\documentclass[12pt,reqno]{amsart}
\usepackage{amsmath,amssymb,amsthm}
\usepackage{mathrsfs}
\usepackage{stmaryrd}

\newtheorem{thm}{Theorem}[section]
\newtheorem{cor}[thm]{Corollary}

\newtheorem{prop}[thm]{Proposition}
\theoremstyle{definition}
\newtheorem{defn}[thm]{Definition}
\newtheorem{nota}[thm]{Notation}
\theoremstyle{remark}
\newtheorem{rem}[thm]{Remark}
\theoremstyle{remark}

\numberwithin{equation}{section}
\date{December 15, 2008.}
\begin{document}
\title[An Empirical Central Limit Theorem in $\mathbf{L}^1$]{An Empirical Central Limit Theorem in $\mathbf{L}^1$ for stationary sequences.}
\author{Sophie Dede}
\address{LPMA, UPMC Universit\'e Paris 6, Case courier 188, 4, Place Jussieu, 75252 Paris Cedex 05, France.}
\email{sophie.dede@upmc.fr} \maketitle

\begin{abstract}
In this paper, we derive asymptotic results for the $\mathbf{L}^1$-Wasserstein distance between the distribution function and the
corresponding empirical distribution function of a stationary sequence. Next, we give some applications to dynamical systems and causal linear processes.
To prove our main result, we give a
Central Limit Theorem for ergodic stationary sequences of random variables with values in $\mathbf{L}^1$.
The conditions obtained are expressed in terms of projective-type conditions. The main tools are martingale approximations.
\end{abstract}
\medskip
\noindent \textit{Mathematics Subject Classifications (2000):} $60$F$17$, $60$G$10$, $62$G$30$. \\
\textit{Key Words:} empirical distribution function, central limit theorem, stationary sequences, Wasserstein distance.
\section{Introduction}
\indent The Kantorovich or $\mathbf{L}^1$-Wasserstein distance between two probability measures $P_1$ and $P_2$ on $\mathbb{R}$
with finite mean, is defined by
\begin{equation*}
d_1(P_1,P_2):=\inf \big \{ \int |x-y| \, d\nu(x,y): \ \nu \in \mathcal{P}(\mathbb{R}^2) \ \mbox{with marginals $P_1$, $P_2$} \big \},
\end{equation*}
where $\mathcal{P}(\mathbb{R}^2)$ is the space of probability measures on $\mathbb{R}^2$. \\
\indent Let $\Lambda_1$ be the space of $1$-Lipschitz functions.
It is well known that $d_1$ can also be written as follows:
\begin{equation*}
d_1(P_1,P_2)
= \int |F_2(t)-F_1(t)| \, dt = \underset{f \in \Lambda_1}{\sup} \Big| \int f dP_1-\int f dP_2 \Big|,
\end{equation*}
where $F_1$ (respectively $F_2$) is the distribution function of $P_1$ (respectively of $P_2$).
\indent Let $(X_i)_{i \in \mathbb{Z}}$ be a stationary sequence of real-valued random variables. In this paper, we are concerned with the Central Limit Theorem (CLT) for the $\mathbf{L}^1$-Wasserstein distance, defined by
\begin{equation}\label{eqnintro}
\int_{\mathbb{R}} |F_n(t)-F_X(t)| \, dt,
\end{equation}
where $F_X$ is the common distribution function of the variables
$X_i$, and $F_n$ is the corresponding empirical distribution function (see Section $3$). \\
\indent In the literature, several previous works on the Kantorovich or $\mathbf{L}^1$-Wasserstein distance, have already been done, for a sequence of i.i.d random variables $X=(X_i)_{i \in \mathbb{Z}}$ ( see for instance del Barrio, Gin\'e and Matr\'an \cite{delbarrioginematran99}).
Recall that if $X$ has the distribution function $F_X$, then the condition
\begin{equation*}
\int_{-\infty}^{\infty} \sqrt{F_X(t)(1-F_X(t))} \, dt < \infty,
\end{equation*}
is equivalent to
\begin{equation*}
\Lambda_{2,1}(X):=\int_0^{\infty} \sqrt{\mathbb{P}(|X|>t)} \, dt<\infty.
\end{equation*}
In their Theorem $2.1$, del Barrio, Gin\'e and Matr\'an \cite{delbarrioginematran99} prove that if $(X_i)_{i \in \mathbb{Z}}$ is i.i.d, then the processes $\sqrt{n}(F_n-F_X)$ converge in law in $\mathbf{L}^1$ to the process $\{B(F(t)),t \in \mathbb{R}\}$, where $B$ is
a Brownian bridge, if and only if $\Lambda_{2,1}(X)<\infty$. Our main result extends Theorem $2.1$ in del Barrio, Gin\'e and Matr\`an \cite{delbarrioginematran99} to the case of stationary sequences, satisfying some appropriate dependence conditions. \\
\indent Before giving the idea of the proof, let us introduce $\mathbf{L}^1(\mu)=\mathbf{L}^1(\mathbf{T},\mu)$, where $\mu$ is a $\sigma$-finite measure,
the Banach space of $\mu$-integrable real functions on $\mathbf{T}$, with the norm $\|.\|_{1,\mu}$, defined by $\|x\|_{1,\mu}=\int_{\mathbf{T}} |x(t)|\, \mu(dt)$.
Let $\mathbf{L}^{\infty}(\mu)$ be its dual space. \\
\indent First, we give the Central Limit Theorem (CLT) for ergodic stationary sequences of martingale differences in $\mathbf{L}^{1}(\mu)$ (see Section 4.1).
Then, by martingale approximation (see for instance Voln\'y \cite{volny93}), we derive a Central Limit Theorem for some ergodic
stationary sequences of $\mathbf{L}^1(\mu)$-valued random variables satisfying some projective criteria. This result allows us to get sufficient conditions to derive
the asymptotic behavior of (\ref{eqnintro}). \\
\indent The paper is organized as follows. In Section 2, we state our main result. In Section 3, we derive the empirical Central Limit Theorem
for statistics of the type (\ref{eqnintro}) for a large class of dependent sequences. In particular, the results apply to unbounded functions of
expanding maps of the interval, and to causal linear processes.

\section{Central Limit Theorem for stationary sequences in $\mathbf{L}^1(\mu)$}
From now, we assume that the ergodic stationary sequence $(X_i)_{i \in \mathbb{Z}}$ of centered random variables with values in $\mathbf{L}^1(\mu)$, is given by $X_i=X_0 \circ \mathbb{T}^i$, where $\mathbb{T}:\Omega \longrightarrow \Omega$ is a bijective
bimeasurable transformation preserving the probability $\mathbb{P}$ on $(\Omega, \mathcal{A})$. Let $S_n=\sum_{j=1}^n X_j$, be the partial sums. For a subfield $\mathcal{F}_0$ satisfying
$\mathcal{F}_0\subseteq \mathbb{T}^{-1}(\mathcal{F}_0)$, let $\mathcal{F}_i=\mathbb{T}^{-i}(\mathcal{F}_0)$.
\begin{nota}
For any integer $p\geq 1$ and for any real random variable $Y$, we denote $\|.\|_p$, the $\mathbf{L}^p$-norm defined by $\|Y\|_p=\mathbb{E}(|Y|^p)^{1/p}$, and
 $\|.\|_{\infty}$ denotes the $\mathbf{L}^{\infty}$-norm, that is the smallest $u$ such that $\mathbb{P}(|Y|>u)=0$.
\end{nota}
Here is our main result:
\begin{thm}\label{thmclt}
Assume that, for any real $t$, $\mathbb{E}(X_0(t)|\mathcal{F}_{-\infty})=0$, $\mathbb{E}(X_0(t)|\mathcal{F}_{\infty})=X_0(t)$ and
\begin{equation}\label{eqnthmclt1}
\int_{\mathbf{T}} \|X_0(t)\|_2 \, \mu(dt)<\infty.
\end{equation}
Let $P_0(X(t))=\mathbb{E}(X(t)|\mathcal{F}_0)-\mathbb{E}(X(t)|\mathcal{F}_{-1})$ and assume that
\begin{equation} \label{eqnthmclt2}
\sum_{k \in \mathbb{Z}} \int_{\mathbf{T}} \|P_0(X_k(t))\|_2 \, \mu(dt)<\infty.
\end{equation}
Then
\begin{equation} \label{eqnthmclt3}
n^{-1/2}\sum_{i=1}^n X_0 \circ \mathbb{T}^i \underset{n \rightarrow \infty}{\longrightarrow} G \ \mbox{in law in $\mathbf{L}^1(\mu)$},
\end{equation}
where $G$ is a $\mathbf{L}^1(\mu)$-valued centered Gaussian random variable
with covariance operator: for any $f \in \mathbf{L}^{\infty}(\mu)$,
\begin{equation} \label{eqnthmclt4}
\Phi_G(f,f)=\mathbb{E}\big (\big(f \big(\sum_{k \in \mathbb{Z}}P_0(X_k) \big) \big)^2 \big)=\sum_{k\in \mathbb{Z}} \mathrm{Cov}(f(X_0),f(X_k)).
\end{equation}
\end{thm}
As a consequence, we have
\begin{cor}\label{corclt1}
Assume that (\ref{eqnthmclt1}) holds.
Moreover, suppose that
\begin{equation} \label{eqncorclt12}
\sum_{n=1}^{\infty} \frac{1}{\sqrt{n}} \int_{\mathbf{T}}\|\mathbb{E}(X_n(t)\mid \mathcal{F}_0)\|_2 \, \mu(dt) < \infty \, ,
\end{equation}
and that
\begin{equation} \label{eqncorclt12bis}
 \sum_{n=1}^{\infty} \frac{1}{\sqrt{n}} \int_{\mathbf{T}}\|X_{-n}-\mathbb{E}(X_{-n}(t) \mid \mathcal{F}_0)\|_2 \, \mu(dt) < \infty \, .
\end{equation}
Then, the conclusion of Theorem \ref{thmclt} holds.
\end{cor}

\section{Applications to the empirical distribution function}
\indent Let $Y=(Y_i)_{i \in \mathbb{Z}}$ be a sequence of real-valued random variables.
We denote
their common distribution function by $F_Y$ and by $F_n$ the corresponding empirical distribution function of $Y$:
\begin{equation*}
\forall \ t \in \mathbb{R}, \ F_n(t)=\frac{1}{n} \sum_{i=1}^n \mathbf{1}_{Y_i \leq t}.
\end{equation*}
\indent Let $\lambda$ be the Lebesgue measure on $\mathbb{R}$. If $\mathbb{E}(|Y_1|)<\infty$, the random variable $X_i(.)=\{t \mapsto \mathbf{1}_{Y_i \leq t}-F_Y(t), t \in \mathbb{R} \}$ may be viewed as
a centered random variable with values in $\mathbf{L}^1(\lambda)$.
\begin{nota}
Let $F_{Y_k \mid \mathcal{F}_0}$ be the conditional function of $Y_k$ given $\mathcal{F}_0$, and let
$F_{Y_k \mid \mathcal{F}_{-1}}$ be the conditional function of $Y_k$ given $\mathcal{F}_{-1}$.
\end{nota}
With these notations, the following equalities are valid:
for every $k$ in $\mathbb{Z}$,
\begin{eqnarray}
\int_{\mathbf{T}} \| P_0(X_k(t))\|_2 \, dt & = & \int_{\mathbf{T}} \|F_{Y_k \mid \mathcal{F}_0}(t)-F_{Y_k \mid \mathcal{F}_{-1}}(t)\|_2 \, dt, \label{eqnappl1} \\
\mbox{and} \ \int_{\mathbf{T}} \|\mathbb{E}(X_k(t)\mid \mathcal{F}_0)\|_2 \, dt & = & \int_{\mathbf{T}} \| F_{Y_k \mid \mathcal{F}_0}(t)-F_{Y}(t)\|_2 \, dt. \label{eqnappl2}
\end{eqnarray}
\subsection{Dependent sequences}
$\newline$
\indent As we shall see in this section, applying Corollary \ref{corclt1}, we can derive sufficient conditions for the convergence in ${\mathbf L}^1(\lambda)$ of the process
$\sqrt n ( F_n - F_Y)$, as soon as the sequence $Y$ satisfies some weak dependence conditions.
Set $\mathcal{F}_0=\sigma(Y_i,i\leq 0)$. We first recall the following dependence coefficients as defined in Dedecker and Prieur \cite{dedeckerprieur05}: for any integer $k \geq 0$,
\begin{equation*}
\tilde{\phi}(k)=\underset{t \in \mathbb{R}}{\sup} \|\mathbb{P}(Y_k \leq t \mid \mathcal{F}_0)-\mathbb{P}(Y_k \leq t)\|_{\infty} \, ,
\end{equation*}
and
\begin{equation*}
\tilde{\alpha}(k)=\underset{t \in \mathbb{R}}{\sup} \|\mathbb{P}(Y_k \leq t \mid \mathcal{F}_0)-\mathbb{P}(Y_k \leq t)\|_1 \, .
\end{equation*}
When the sequence $Y$ is $\tilde{\phi}$-dependent, the following result holds:
\begin{prop} \label{propphimelange}
Assume that
\begin{equation}\label{eqnpropphimelange1}
\sum_{k \geq 1} \sqrt{\frac{\tilde{\phi}(k)}{k}} < \infty \ \mbox{and} \ \int_{0}^{\infty} \sqrt{\mathbb{P}(|Y|>t)} \, dt <\infty,
\end{equation}
then $\{t\mapsto \sqrt{n}(F_n(t)-F_Y(t)), t \in \mathbb{R}\}$ converges in $\mathbf{L}^1(\lambda)$, to
a centered Gaussian random variable, with covariance function:
for any $f, \ g \in \mathbf{L}^{\infty}(\lambda)$,
\begin{equation}\label{eqnpropphimelange10}
\Phi_{\lambda}(f,g) = \int_{\mathbb{R}^2} f(s)g(t) C(s,t) \, dt \, ds
\end{equation}
with
\begin{equation*}
C(s,t) = F_Y(t \wedge s)-F_Y(t) F_Y(s) + 2 \sum_{k \geq 1} ( \mathbb{P}(Y_0 \leq t, Y_k \leq s )-F_Y(t)F_Y(s)).
\end{equation*}
\end{prop}
\begin{rem} Proposition \ref{propphimelange} is also true with the $\phi$-mixing coefficient of Ibragimov \cite{ibragimov62}.
Notice that this result contains the i.i.d case, developed in del Barrio, Gin\'e and Matr\'an \cite{delbarrioginematran99}.
\end{rem}
\bigskip

Before giving sufficient conditions when the sequence $Y$ is $\tilde{\alpha}$-dependent, we first recall the following definition:
\begin{defn}
For any nonnegative and integrable random variable $Y$, define the quantile function $Q_Y$ of $|Y|$, that is the cadlag inverse of the tail function
$x \rightarrow \mathbb{P}(|Y|> x)$.
\end{defn}
\begin{prop} \label{propalphamelange}
Assume that
\begin{equation} \label{eqnpropalphamelange1}
\sum_{k \geq 1} \frac{1}{\sqrt{k}} \int_0^{\tilde{\alpha}(k)} \frac{Q_Y(u)}{\sqrt{u}} \, du < \infty,
\end{equation}
then the conclusion of Proposition \ref{propphimelange} holds.
\end{prop}
\begin{rem} Notice that Proposition \ref{propalphamelange} is also true with strong $\alpha$-mixing coefficients of Rosenblatt \cite{rosenblatt56}. Notice also that (\ref{eqnpropalphamelange1}) is equivalent to
\end{rem}
\begin{equation}\label{eqnpropalphamelange2}
\sum_{k \geq 1} \frac{1}{\sqrt{k}} \int_0^{\infty} \sqrt{\tilde{\alpha}(k)}\wedge \sqrt{\mathbb{P}(|Y|>t)} \, dt < \infty.
\end{equation}
\subsubsection{Application to expanding maps.}
\indent Let $T$ be a map from $[0,1]$ to $[0,1]$ preserving a probability $\mu$ on $[0,1]$. Recall that  the Perron-Frobenius operator $K$ from $\mathbf{L}^1(\mu)$ to $\mathbf{L}^{1}(\mu)$ is defined via the equality: for any $h \in \mathbf{L}^1(\mu)$ and $f \in \mathbf{L}^{\infty}(\mu)$,
\begin{equation*}
\int_0^1(Kh)(x) f(x) \mu(dx)=\int_0^1 h(x) (f \circ T)(x) \mu(dx).
\end{equation*}
Here we are interested by giving sufficient conditions for the convergence in ${\mathbf L}^1(\lambda)$ of the empirical distribution function associated to $F_Y$ where the random variables $(Y_i)_{i\in {\mathbb Z}}$ are defined as follows: for a given monotonic function $f$, let
\begin{equation} \label{defYk}
Y_k=f\circ T^k \, .
\end{equation}
In fact since on the probability $([0,1],\mu)$, the random variable $(T,T^2, ..., T^n)$ is distributed as $(Z_n,Z_{n-1},...,Z_1)$, where $(Z_i)_{i \geq 0}$ is a stationary Markov chain with invariant measure $\mu$ and transition Kernel $K$ (see Lemma $XI.3$ in Hennion and Herv\'e \cite{hennionherve01}), the convergence
in ${\mathbf L}^1(\lambda)$ of the empirical distribution function associated to $F_Y$ is reduced to the one of the empirical distribution function associated to $F_{f(Z)}$.

\medskip

In this section we consider two cases: first the case of a class of BV-contracting maps and secondly the case of a class of intermittent maps.

\medskip

\noindent\textit{a) The case of BV-contracting maps.}
Let $BV$ be the class of bounded variation functions from $[0,1]$ to $\mathbb{R}$.
For any $h \in BV$, denote by
$\|dh\|$ the variation norm of the measure $dh$.
A Markov kernel $K$ is said to be $BV$-contracting if there exist $C >0$ and $\rho \in [0,1[ $ such that
\begin{equation}\label{eqn110000}
\| dK^n(h)\| \leq C \rho^n \|dh\|.
\end{equation}
A map $T$ is then said to be BV-contracting if its Perron-Frobenius operator $K$ is BV-contracting (see for instance Dedecker and Prieur \cite{dedeckerprieur05}, for more details and examples of maps which are BV-contracting).

In this case, the following result holds:

\begin{cor}\label{exphi}
If $T$ is BV-contracting and $f$ is a monotonic function from $]0,1[$ to $\mathbb{R}$ satisfying $\int_0^{\infty} \sqrt{\lambda(|f|>t)} \, dt <\infty$, then the conclusion of Proposition \ref{propphimelange} holds for the sequence $(Y_k)_{k \in {\mathbb Z}}$ where $Y_k$ is defined by (\ref{defYk}).
\end{cor}
\begin{rem}In the particular case when $f$ is positive and non increasing on $]0,1[$, with $f(x) \leq D x^{-a}$ for some $a>0$ and $D$ a constant, we get that
\begin{equation*}
\int_0^{\infty} \sqrt{\lambda(|f|>t)} \, dt \leq C_2 \int_1^{\infty} \frac{1}{t^{1/(2a)}} \, dt,
\end{equation*}
where $C_2$ is a constant. Consequently, Corollary \ref{exphi} holds as soon as $a<\frac{1}{2}$ holds.
\end{rem}

\medskip

\noindent \textit{b) Application to intermittent maps.}  For $\gamma$ in $]0,1[$, we consider the intermittent map $T_{\gamma}$ from $[0,1]$ to $[0,1]$, studied for instance by Liverani, Saussol and Vaienti \cite{liveranisaussolvaienti99}, which is a
modification of the Pomeau-Manneville map \cite{pomeaumanneville80}:
\begin{displaymath}
T_{\gamma}=\left\{ \begin{array}{cc}
          x (1+2^{\gamma} x^{\gamma}) & \mbox{if $x \in [0,1/2[$} \\
          2x-1 & \mbox{if $x \in [1/2,1]$.}
          \end{array}
          \right.
\end{displaymath}
We denote by $\nu_{\gamma}$ the unique $T_{\gamma}$-probability measure on $[0,1]$ and by $K_{\gamma}$ the Perron-Frobenius operator of $T_{\gamma}$
with respect to $\nu_{\gamma}$. For these maps, we obtain the following result:

\begin{cor}\label{coralpha}
For $\gamma$ in $]0,1[$, if $T_{\gamma}$ is an intermittent map and $f$ is a monotonic function from $]0,1[$ to $\mathbb{R}$, satisfying
\begin{equation}\label{eqnalphamixing1}
\sum_{k\geq 1} \frac{1}{\sqrt{k}} \int_{0}^{\infty} \frac{1}{k^{\frac{1-\gamma}{2\gamma}}} \wedge \sqrt{\nu_{\gamma}(|f|>t)}\, dt <\infty,
\end{equation}
then the conclusion of Proposition \ref{propalphamelange} holds for the sequence $(Y_k)_{k \in {\mathbb Z}}$ where $Y_k$ is defined by (\ref{defYk}).
\end{cor}
\begin{rem} \label{remG}
In the particular case when $f$ is positive and non increasing on $]0,1[$, with $f(x) \leq D x^{-a}$ for some $a>0$ and $g_{\nu_{\gamma}}$ the density of $\nu_{\gamma}$ such that $g_{\nu_{\gamma}}(x) \leq V_{\gamma} x^{-\gamma}$ where $V(\gamma)$ is a constant, we can prove that
(\ref{eqnalphamixing1}) holds
as soon as $a < \frac{1}{2}-\gamma$ does (see Section
\ref{proofremG}). In his comment after Theorem $3$, Gou\"ezel \cite{gouezel04} proved that if $f(x)=x^{-a}$, then $n^{-1/2} \sum_{k=1}^n (f \circ T_{\gamma}^i-\nu_{\gamma}(f))$
converges to a normal law if $a<1/2-\gamma$, and that there is a convergence to a stable law (with a different normalization) if $a>1/2-\gamma$.
This example shows that our condition is close to optimality.
\end{rem}

\subsection{Causal linear processes}
$\newline$
\indent We  focus here on the stationary sequence
\begin{equation}\label{defcausallinearadapt}
Y_k=\sum_{j \geq 0} a_j \varepsilon_{k-j},
\end{equation}
where $(\varepsilon_i)_{i \in \mathbb{Z}}$ is a sequence of real-valued i.i.d random variables in $\mathbf{L}^2$ and $\sum_{j \geq 0} |a_j|<+\infty$.
\begin{cor}\label{propcausallinearprocessadapt}
Assume that, $\varepsilon_0$ has a density bounded by $K$ and that $|a_0|\not=0$.
Moreover, assume that
\begin{equation}\label{eqnpropcausallinearprocessadapt1}
\sum_{k \geq 0} \int_0^{(a_k)^2} \frac{Q_{\mid Y_0 \mid}(u)}{\sqrt{u}} \, du<\infty.
\end{equation}
Then the conclusion of Proposition \ref{propphimelange} holds for the sequence $(Y_k)_{k \in {\mathbb Z}}$ where $Y_k$ is defined by (\ref{defcausallinearadapt}).
\end{cor}
\begin{rem}\label{remRio}
Since $\sum_{j\geq 0}|a_j|<\infty$, (\ref{eqnpropcausallinearprocessadapt1}) is true provided that
\begin{equation}\label{eqncausallinearprocess3}
\sum_{k\in \mathbb{Z}} \int_0^{\mid a_k \mid^2} \frac{Q_{\mid \varepsilon_0 \mid }(u)}{\sqrt{u}} \, du < \infty.
\end{equation}
\end{rem}
As a consequence, we get the following result
\begin{cor} \label{propmom}Assume either Item 1 or 2 below:
\begin{enumerate}
\item[$1.$] for some $r>2$, the i.i.d random variables $(\varepsilon_i)_{i\in \mathbb{Z}}$ are in $\mathbf{L}^r$ and  $\varepsilon_0$ has a density bounded by $K$. In addition  $|a_0|\not=0$ and
\begin{equation}\label{eqnmom}
\sum_{k \geq 0}  k^{1/(r-1)}|a_k|^{(r-2)/(r-1)} < \infty.
\end{equation} \\
\item[$2.$] for some $r>2$,
\begin{equation*}
\forall \ x >0, \ \mathbb{P}(|\varepsilon|>x)\leq \big( \frac{c}{x} \big)^r \mbox{ where $c$ is a positive constant} \, ,
\end{equation*}
and $\varepsilon_0$ has a density bounded by $K$. In addition  $|a_0|\not=0$ and
\begin{equation}\label{eqnqueue}
\sum_{k \geq 0} |a_k|^{1-2/r}<\infty.
\end{equation}
\end{enumerate}
Then the conclusion of Proposition \ref{propphimelange} holds for the sequence $(Y_k)_{k \in {\mathbb Z}}$ where $Y_k$ is defined by (\ref{defcausallinearadapt}).
\end{cor}

\section{Proofs}
\subsection{Central Limit Theorem for $\mathbf{L}^1(\mu)$-valued martingale differences}
$\newline$
\indent We extend to $\mathbf{L}^1(\mu)$-valued martingale differences, a result of Jain \cite{jain77}, which is,
for a sequence of i.i.d centered $\mathbf{L}^1(\mu)$-valued random variables $X=(X_i)_{i\in \mathbb{Z}}$,
the Central Limit Theorem holds if and only if
$\int_{\mathbf{T}} (\mathbb{E}(X_1(t)^2))^{1/2} \, \mu(dt) <\infty$.\\
\begin{thm} \label{thmcltdiffmart}
Let $(M_i)_{i\in \mathbb{Z}}$ be a sequence of stationary ergodic martingale differences
with values in $\mathbf{L}^1(\mu)$ such that $M_i=M_0 \circ \mathbb{T}^i$.
Assume that
\begin{equation} \label{eqnthmcltdiffmart1}
\int_{\mathbf{T}} \|M_0(t)\|_2 \, \mu(dt)<\infty.
\end{equation}
Then
\begin{equation} \label{eqnthmcltdiffmart2}
n^{-1/2} \sum_{i=1}^n M_0 \circ \mathbb{T}^i \underset{n \rightarrow \infty}{\longrightarrow} G \ \mbox{in law},
\end{equation}
where $G$ is a $\mathbf{L}^1(\mu)$-valued centered Gaussian random variable with covariance function: for any $f$ in $\mathbf{L}^{\infty}(\mu)$, $\Phi_G(f,f)=\mathbb{E}(f^2(M_0))$.
\end{thm}
\noindent \textit{Proof of Theorem \ref{thmcltdiffmart}}.
Before giving the proof, we recall some notations and definitions used in the proof.
\begin{nota}Using the notations of Jain \cite{jain77}, we consider for any real separable Banach space $\mathbb{B}$ and its dual $\mathbb{B}'$,
\begin{equation*}
WM_0^2=\{ \nu \ \mbox{ probability measures on $\mathbb{B}$:} \ \int |f|^2 \, d\nu < \infty, \ \int f \, d\nu =0, \ \forall \ f \in \mathbb{B}'\}.
\end{equation*}
\end{nota}
\begin{defn}
If $\nu \in WM_0^2$, its covariance kernel $\Phi_{\nu}$ is given by: for any $f, g \in \mathbb{B}'$,
\begin{equation*}
\Phi_{\nu}(f,g)=\int fg \, d\nu.
\end{equation*}
\end{defn}
\begin{defn}
$\mu \in WM_0^2$ is pregaussian if there is a Gaussian measure $\nu$, such that $\Phi_{\nu}=\Phi_{\mu}$.
\end{defn}
By the classical Linderberg's theorem for stationary ergodic
martingale differences in Billingsley \cite{billingsley61}, $n^{-1/2} f(S_n)$ converges in law to the centered gaussian random variable $Z$ in
$\mathbb{R}$, with variance $\mathbb{E}(f^2(M_0))$, for each $f\in \mathbf{L}^{\infty}(\mu)$.
Now, we have to prove that the distribution of $n^{-1/2}S_n$ is relatively compact. As $\mathbf{L}^1(\mu)$ is of cotype $2$, we use the same approach
as in the proof of Theorem $6.4$ in de Acosta, Araujo and Gin\'e \cite{deacostaaraujogine78}.\\
By stationarity, it follows that
\begin{equation*}
\mathbb{E}\Big(\Big[f\Big(\frac{S_n}{\sqrt{n}}\Big)\Big]^2\Big)=\frac{1}{n} \sum_{i=1}^n \mathbb{E}([f(M_i)]^2)=\mathbb{E}([f(M_0)]^2).
\end{equation*}
By Theorem $11$ in Jain \cite{jain77}, if $(\ref{eqnthmcltdiffmart1})$ holds, then $M_0$ is pregaussian, so there exists a $\mathbf{L}^1(\mu)$-valued centered Gaussian
random variable $G$ ( or a Gaussian measure $\gamma$ on $\mathbf{L}^1(\mu)$), with covariance operator $\Phi_G$, such that, for any $f$ in $\mathbf{L}^{\infty}(\mu)$,
\begin{equation*}
\mathbb{E}([f(G)]^2)=\Phi_G(f,f)=\mathbb{E}([f(M_0)]^2).
\end{equation*}
By Theorem $5.6$ in de Acosta, Araujo and Gin\'e \cite{deacostaaraujogine78}, every centered Gaussian measure on $\mathbf{L}^1(\mu)$,
is strongly Gaussian which means that there exist a Hilbert space $\mathbb{H}$, a continuous linear map $\mathbf{M}: \mathbb{H} \rightarrow \mathbf{L}^1(\mu)$ and
a tight centered Gaussian measure $\nu$ on $\mathbb{H}$ such that $\gamma=\nu \circ \mathbf{M}^{-1}$.
Therefore, we can apply Theorem $6.2$ \cite{deacostaaraujogine78}, with $\mathfrak{K}=\{ \xi \ \mbox{a probability measure on $\mathbf{L}^1 (\mu)$ such that} \ \Phi_{\xi}(f,f) \leq \Phi_{\gamma}(f,f), \mbox{for all $f \in \mathbf{L}^{\infty}(\mu)$}\}$, so $\mathfrak{K}$ is relatively compact. We have proved that the distribution of $n^{-1/2}S_n$ is relatively compact.
\hfill
$\square$
\subsection{Proof of Theorem \ref{thmclt}}
$\newline$
\indent We construct the martingale
\begin{equation*}
M_n=\sum_{i=1}^n M_0 \circ \mathbb{T}^i,
\end{equation*}
where $M_0=\sum_{k \in \mathbb{Z}} P_0 (X_k)$. Notice that $(M_0 \circ \mathbb{T}^i)_{i \in \mathbb{Z}}$ is a sequence of a stationary ergodic martingale differences. By triangle inequality,
\begin{eqnarray*}
\int_{\mathbf{T}} \|M_0(t)\|_2 \, \mu(dt) & = & \int_{\mathbf{T}} \big \|\sum_{k \in \mathbb{Z}} P_0(X_k(t)) \big \|_2 \, \mu(dt) \\
                                   & \leq & \sum_{k \in \mathbb{Z}} \int_{\mathbf{T}} \|P_0(X_k(t))\|_2 \, \mu(dt) < \infty.
\end{eqnarray*}
Applying Theorem \ref{thmcltdiffmart}, we infer that
\begin{equation*}
\frac{1}{\sqrt{n}} M_n \underset{n\rightarrow \infty}{\longrightarrow} G \ \mbox{in law in $\mathbf{L}^1(\mu)$},
\end{equation*}
where $G$ is a $\mathbf{L}^1(\mu)$-valued centered Gaussian random variable such that $\Phi_G(f,f)=\mathbb{E}([f(M_0)]^2)$, for any $f \in \mathbf{L}^{\infty}(\mu)$.\\
\indent To conclude the proof, it suffices to prove that
\begin{equation}\label{eqndemthmclt1}
\underset{n \longrightarrow \infty}{\lim}  \int_{\mathbf{T}} \Big \| \frac{S_n(t)}{\sqrt{n}} -\frac{1}{\sqrt{n}} \sum_{i=1}^n M_0(t) \circ \mathbb{T}^i \Big \|_2 \, \mu(dt)=0 .
\end{equation}
The proof is inspired by the proof of Theorem $1$ in Dedecker, Merlev\`ede and Voln\'y \cite{dedeckermerlevedevolny07}.
By triangle inequality,
\begin{align}
& \int_{\mathbf{T}} \Big \|  \frac{S_n(t)}{\sqrt{n}} - \frac{1}{\sqrt{n}} \sum_{i=1}^n M_0(t)\circ \mathbb{T}^i \Big \| _2 \, \mu(dt) \notag   \\
=& \int_{\mathbf{T}} \Big \| \frac{S_n(t)}{\sqrt{n}} - \frac{\mathbb{E}(S_n(t) \mid \mathcal{F}_n)}{\sqrt{n}} + \frac{\mathbb{E}(S_n(t) \mid \mathcal{F}_n)}{\sqrt{n}} + \frac{\mathbb{E}(S_n(t) \mid \mathcal{F}_0)}{\sqrt{n}} \notag  \\
& \qquad \qquad \qquad \qquad \qquad - \frac{\mathbb{E}(S_n(t) \mid \mathcal{F}_0)}{\sqrt{n}} - \frac{1}{\sqrt{n}} \sum_{i=1}^n M_0(t)\circ \mathbb{T}^i \Big \| _2 \, \mu(dt) \notag \\
 \leq & \int_{\mathbf{T}}\Big \| \frac{S_n(t)}{\sqrt{n}} - \frac{\mathbb{E}(S_n(t) \mid \mathcal{F}_n)}{\sqrt{n}} \Big \| _2 \, \mu(dt) \notag \\
& \qquad + \int_{\mathbf{T}} \Big \| \frac{\mathbb{E}(S_n(t) \mid \mathcal{F}_n)}{\sqrt{n}} - \frac{\mathbb{E}(S_n(t) \mid \mathcal{F}_0)}{\sqrt{n}} - \frac{1}{\sqrt{n}} \sum_{i=1}^n M_0(t)\circ \mathbb{T}^i \Big \| _2 \, \mu(dt) \notag\\
 & \qquad \qquad \qquad \qquad \qquad+ \int_{\mathbf{T}} \Big  \| \frac{\mathbb{E}(S_n(t) \mid \mathcal{F}_0)}{\sqrt{n}} \Big \| _2 \, \mu(dt) \label{eqndemthmclt11}.
\end{align}
It suffices to prove that each term of the right-hand side in Inequality $(\ref{eqndemthmclt11})$ tends to $0$, as $n$ tends to infinity.
Let us first control the second term. Since
\begin{equation*}
\mathbb{E}(S_n(t) \mid \mathcal{F}_n)- \mathbb{E}(S_n(t) \mid \mathcal{F}_0)=\sum_{i=1}^n \sum_{k=1}^n P_i(X_k(t)),
\end{equation*}
it follows that, by stationarity and orthogonality,
\begin{align}
&\int_{\mathbf{T}} \Big \| \frac{\mathbb{E}(S_n(t) \mid \mathcal{F}_n)- \mathbb{E}(S_n(t) \mid \mathcal{F}_0)}{\sqrt{n}} - \frac{1}{\sqrt{n}} \sum_{i=1}^n M_0(t)\circ \mathbb{T}^i \Big \| _2 \, \mu(dt) \notag  \\
=&  \int_{\mathbf{T}} \sqrt{ \Big \| \frac{\mathbb{E}(S_n(t) \mid \mathcal{F}_n)-\mathbb{E}(S_n(t) \mid \mathcal{F}_0)}{\sqrt{n}} - \frac{1}{\sqrt{n}} \sum_{i=1}^n M_0(t)\circ \mathbb{T}^i \Big \| _2^2 } \, \mu(dt) \notag  \\
=&  \int_{\mathbf{T}} \sqrt{\frac{1}{n} \sum_{i=1}^n \big \| \sum_{k=1}^n P_0(X_{k-i}(t))-M_0(t) \big \| _2^2} \, \mu(dt) \notag \\
=& \ \int_{\mathbf{T}} \sqrt{\frac{1}{n} \sum_{i=1}^n \big \| \sum_{j=1-i}^{n-i} P_0(X_j(t))-M_0(t) \big \| _2^2 } \, \mu(dt) \notag  \\
\leq &  \sqrt{2} \Big [  \int_{\mathbf{T}} \sqrt{\frac{1}{n} \sum_{i=1}^n \big \| \sum_{j\leq -i}P_0(X_j(t))\big \| _2^2} \, \mu(dt) +  \int_{\mathbf{T}} \sqrt{\frac{1}{n} \sum_{i=1}^n \big \| \sum_{j \geq n-i+1} P_0(X_j(t)) \big \| _2^2} \, \mu(dt) \Big]. \label{eqndemthmclt2}
\end{align}
Splitting the sum on $i$ of the first term in the right-hand side of Inequality $(\ref{eqndemthmclt2})$, we get that
\begin{align*}
 & \int_{\mathbf{T}} \sqrt{\frac{1}{n} \sum_{i=1}^n \big \| \sum_{j \leq -i} P_0(X_j(t))\big  \| ^2_2} \, \mu(dt) \\
=& \int_{\mathbf{T}} \sqrt{\frac{1}{n} \sum_{i=1}^N \big \| \sum_{j \leq -i} P_0(X_j(t))\big \| ^2_2+ \frac{1}{n} \sum_{i=N+1}^n \big \| \sum_{j \leq -i} P_0(X_j(t))\big \| ^2_2 } \, \mu(dt)  \\
\leq & \int_{\mathbf{T}} \sqrt{\frac{1}{n} \sum_{i=1}^N \big \| \sum_{j \leq -i} P_0(X_j(t))\big \| _2^2} \, \mu(dt) + \int_{\mathbf{T}} \sqrt{\frac{1}{n} \sum_{i=N+1}^n \big \| \sum_{j \leq -i} P_0(X_j(t))\big \| ^2_2} \, \mu(dt).
\end{align*}
Fubini entails that
\begin{eqnarray*}
\int_{\mathbf{T}} \sqrt{\frac{1}{n} \sum_{i=1}^N \big \| \sum_{j \leq -i} P_0(X_j(t))\big \| ^2_2} \, \mu(dt) & \leq & \frac{1}{\sqrt{n}} \int_{\mathbf{T}} \sum_{i=1}^N \big \| \sum_{j \leq -i} P_0(X_j(t))\big \| _2 \, \mu(dt) \\
                                                                                  & \leq & \frac{1}{\sqrt{n}} \int_{\mathbf{T}} N  \sum_{j \leq -1} \big \|P_0(X_j(t))\big\| _2 \, \mu(dt) \\
                                                                                  & \leq & \frac{1}{\sqrt{n}}  N \sum_{j \in \mathbb{Z}} \int_{\mathbf{T}} \| P_0(X_j(t))\| _2 \, \mu(dt) \underset{n\rightarrow \infty}{\longrightarrow} 0.
\end{eqnarray*}
Moreover, since
\begin{eqnarray*}
\int_{\mathbf{T}} \sqrt{\frac{1}{n} \sum_{i=N+1}^n \big \| \sum_{j \leq -i} P_0(X_j(t)) \big \| ^2_2} \, \mu(dt)
 &\leq & \int_{\mathbf{T}} \sqrt{\frac{1}{n}(n-N) \big( \sum_{j \leq -N} \| P_0(X_j(t))\| _2 \big)^2} \, \mu(dt) \\
&\leq & \int_{\mathbf{T}} ( \sum_{j \leq -N} \| P_0(X_j(t))\| _2) \, \mu(dt)  \\
&\leq & \sum_{j \leq -N} \int_{\mathbf{T}} \| P_0(X_j(t))\| _2 \, \mu(dt), \\
\end{eqnarray*}
we infer by $(\ref{eqnthmclt2})$ that
\begin{equation*}
\underset{N \rightarrow \infty}{\lim} \, \underset{n \rightarrow \infty}{\limsup} \, \int_{\mathbf{T}} \sqrt{\frac{1}{n}\sum_{i=N+1}^n \big \| \sum_{j \leq -i} P_0(X_j(t))\|_2^2} \, \mu(dt)=0 \, .
\end{equation*}
Whence
\begin{equation*}
\int_{\mathbf{T}} \sqrt{\frac{1}{n} \sum_{i=1}^n \big \| \sum_{j \leq -i} P_0(X_j(t))\big \| _2^2} \, \mu(dt) \underset{n \rightarrow \infty}{\longrightarrow } 0.
\end{equation*}
\noindent In the same way, splitting the sum on $i$ of the second term in the right-hand side of Inequality $(\ref{eqndemthmclt2})$, we derive that
\begin{align*}
& \int_{\mathbf{T}} \sqrt{\frac{1}{n} \sum_{i=1}^n \big \| \sum_{j \geq n-i+1} P_0(X_j(t)) \big \| _2^2} \, \mu(dt) \\
= & \int_{\mathbf{T}} \sqrt{ \frac{1}{n} \sum_{i=1}^{n-N} \big \| \sum_{j\geq n-i+1} P_0(X_j(t))\big \| _2^2+ \frac{1}{n} \sum_{i=n-N+1}^n \big \| \sum_{j \geq n-i+1} P_0(X_j(t))\big \| _2^2 } \, \mu(dt) \\
\leq &  \int_{\mathbf{T}} \sqrt{\frac{1}{n} \sum_{i=1}^{n-N} \big \| \sum_{j \geq n-i+1} P_0(X_j(t))\big \| _2^2} \, \mu(dt) + \int_{\mathbf{T}} \sqrt{\frac{1}{n} \sum_{i=n-N+1}^n \big \| \sum_{j \geq n-i+1} P_0(X_j(t))\big \| _2^2} \, \mu(dt).
\end{align*}
Since
\begin{eqnarray*}
\int_{\mathbf{T}} \sqrt{\frac{1}{n} \sum_{i=1}^{n-N} \big \| \sum_{j \geq n-i+1} P_0(X_j(t)) \big \| ^2_2} \, \mu(dt) & \leq & \int_{\mathbf{T}} \sqrt{\frac{(n-N)}{n}} \Big (\sum_{j \geq N+1} \| P_0(X_j(t))\| _2 \Big) \, \mu(dt)  \\
                                                                                         & \leq & \sum_{j \geq N+1} \int_{\mathbf{T}} \| P_0(X_j(t))\| _2 \, \mu(dt),  \\
\end{eqnarray*}
and
\begin{eqnarray*}
\int_{\mathbf{T}} \sqrt{\frac{1}{n} \sum_{i=n-N+1}^n \big \| \sum_{j \geq n-i+1} P_0(X_j(t)) \big \| _2^2} \, \mu(dt) & \leq & \int_{\mathbf{T}} \frac{1}{\sqrt{n}} \sum_{i=n-N+1}^n \big \| \sum_{j \geq n-i+1} P_0(X_j(t)) \big \| _2 \, \mu(dt)  \\
                                                                                        & \leq & \frac{N}{\sqrt{n}} \int_{\mathbf{T}} \sum_{j \in \mathbb{Z}} \| P_0(X_j(t))\| _2 \, \mu(dt),
\end{eqnarray*}
we deduce by $(\ref{eqnthmclt2})$, that
\begin{equation*}
\underset{n \rightarrow \infty}{\lim} \int_{\mathbf{T}} \sqrt{\frac{1}{n} \sum_{i=1}^n  \big \| \sum_{j \geq n-i+1} P_0(X_j(t)) \big \| ^2_2} \, \mu(dt) =0.
\end{equation*}
Consequently, we derive that
\begin{equation*}
\int_{\mathbf{T}} \Big \| \mathbb{E}\Big(\frac{S_n(t)}{\sqrt{n}} \Big | \mathcal{F}_n \Big)-\mathbb{E}\Big(\frac{S_n(t)}{\sqrt{n}}\Big | \mathcal{F}_0\Big)-\frac{1}{\sqrt{n}}\sum_{i=1}^n M_0(t)\circ\mathbb{T}^i \Big\|_2 \, \mu(dt) \underset{n\rightarrow \infty}{\longrightarrow} 0.
\end{equation*}
To prove that the last term of Inequality $(\ref{eqndemthmclt11})$ tends to $0$ as $n$ tends to infinity, we first write that
\begin{eqnarray*}
 \int_{\mathbf{T}} \Big \| \mathbb{E}\Big (\frac{S_n(t)}{\sqrt{n}} \Big | \mathcal{F}_0 \Big ) \Big \| _2 \, \mu(dt)
 &\leq & \frac{1}{\sqrt{n}} \int_{\mathbf{T}} \big \| \sum_{k=1}^N \mathbb{E}(X_k(t) \mid \mathcal{F}_0) \big \| _2 \, \mu(dt)  \\
 & & \qquad + \frac{1}{\sqrt{n}} \int_{\mathbf{T}} \big \| \sum_{k=N+1}^n \mathbb{E}(X_k(t) \mid \mathcal{F}_0 ) \big \| _2 \, \mu(dt).
\end{eqnarray*}
By orthogonality,
\begin{eqnarray*}
\big \| \sum_{k=N+1}^n \mathbb{E}(X_k(t) \mid \mathcal{F}_0 ) \big \| _2^2 & =& \sum_{k=N+1}^n \sum_{l=N+1}^n \mathbb{E}( \mathbb{E}(X_k(t) \mid \mathcal{F}_0) \mathbb{E}(X_l(t) \mid \mathcal{F}_0)) \\ \\
                                                                 & = & \sum_{k=N+1}^n \sum_{l=N+1}^n \mathbb{E}( \sum_{m=0}^ {\infty} P_{-m}(X_k(t)) P_{-m}(X_l(t))). \\
\end{eqnarray*}
Using Cauchy-Schwarz inequality and stationarity, it follows that
\begin{eqnarray*}
 \frac{1}{n} \big \| \sum_{k=N+1}^n \mathbb{E}(X_k(t) \mid \mathcal{F}_0) \big \|  _2^2
 &\leq & \frac{1}{n} \sum_{m=0}^{\infty} \sum_{k=N+m+1}^{n+m} \sum_{l=N+1+m}^{n+m} \| P_0(X_k(t))\| _2 \| P_0(X_l(t))\| _2 \\
  &\leq & \Big( \sum_{k=N+1}^{\infty} \| P_0(X_k(t))\| _2 \Big)^2 \, .
\end{eqnarray*}
Consequently
\begin{equation*}
\frac{1}{\sqrt{n}} \int_{\mathbf{T}} \big \| \sum_{k=N+1}^n \mathbb{E}(X_k(t) \mid \mathcal{F}_0) \big \| _2 \, \mu(dt) \leq \sum_{k=N+1}^{\infty} \int_{\mathbf{T}} \| P_0(X_k(t))\| _2 \, \mu(dt) \, ,
\end{equation*}
and by $(\ref{eqnthmclt2})$, it follows that
\begin{equation} \label{eqnthmclt111}
\underset{N\rightarrow \infty}{\lim} \underset{n\rightarrow\infty}{\limsup} \ \frac{1}{\sqrt{n}} \int_{\mathbf{T}} \big \|\sum_{k=N+1}^n \mathbb{E}(X_k(t)\mid \mathcal{F}_0)\big\|_2 \, \mu(dt)=0.
\end{equation}
On the other hand, by stationarity,
\begin{eqnarray*}
\big \| \sum_{k=1}^N \mathbb{E}(X_k(t) \mid \mathcal{F}_0) \big \| _2 & = & \big \| \sum_{k=1}^N \sum_{i \in \mathbb{Z}} \mathbb{E}(P_i(X_k(t)) \mid \mathcal{F}_0) \big \| _2  \\
                                                            & \leq & \sum_{k=1}^N \big ( \sum_{i \in \mathbb{Z}} \| P_i(X_k(t))\| _2 \big )  \\
                                                            & \leq & N \sum_{i \in \mathbb{Z}} \| P_0(X_i(t))\|_2 \, .
\end{eqnarray*}
Hence by $(\ref{eqnthmclt2})$, we get that
\begin{equation}\label{eqnthmclt112}
\lim_{n \rightarrow \infty}\int_{\mathbf{T}} \frac{1}{\sqrt{n}} \big \| \sum_{k=1}^{N} \mathbb{E}(X_k(t) \mid \mathcal{F}_0) \big \| _2 \, \mu(dt)  =0 \, .
\end{equation}
Therefore, $(\ref{eqnthmclt111})$ and $(\ref{eqnthmclt112})$ imply that
\begin{equation*}
\underset{n \rightarrow \infty}{\lim} \ \int_{\mathbf{T}} \Big \|\mathbb{E}\Big(\frac{S_n(t)}{\sqrt{n}} \Big | \mathcal{F}_0 \Big) \Big\|_2 \, \mu(dt)=0.
\end{equation*}
To prove that the first term of Inequality $(\ref{eqndemthmclt11})$ tends to $0$ as $n$ tends to infinity, we write that
\begin{eqnarray*}
\int_{\mathbf{T}} \Big \| \frac{S_n(t)}{\sqrt{n}}- \mathbb{E}\Big( \frac{S_n(t)}{\sqrt{n}} \Big | \mathcal{F}_n\Big)\Big \| _2 \, \mu(dt)
& \leq & \frac{1}{\sqrt{n}} \int_{\mathbf{T}} \Big \| \sum_{k=1}^{n-N} [X_k(t)-\mathbb{E}(X_k(t) \mid \mathcal{F}_n)] \Big \| _2 \, \mu(dt) \\
& & \quad + \frac{1}{\sqrt{n}} \int_{\mathbf{T}} \Big \| \sum_{k=n-N+1}^n [X_k(t)-\mathbb{E}(X_k(t) \mid \mathcal{F}_n)] \Big \| _2 \, \mu(dt).
 \end{eqnarray*}
By orthogonality,
\begin{align*}
&\big \| \sum_{k=1}^{n-N} [X_k(t)-\mathbb{E}(X_k(t) \mid \mathcal{F}_n)]\big \| ^2_2  \\
=&\sum_{k=1}^{n-N} \sum_{l=1}^{n-N} \mathbb{E}([X_k(t)-\mathbb{E}(X_k(t) \mid \mathcal{F}_n)][X_l(t)-\mathbb{E}(X_l(t) \mid \mathcal{F}_n)]) \\
=&\sum_{k=1}^{n-N} \sum_{l=1}^{n-N} [\mathbb{E}(X_k(t)X_l(t))-\mathbb{E}(\mathbb{E}(X_k(t) \mid \mathcal{F}_n)\mathbb{E}(X_l(t) \mid \mathcal{F}_n))] \\
=&  \sum_{k=1}^{n-N} \sum_{l=1}^{n-N} \mathbb{E}\big (\sum_{m=n+1}^{\infty} P_m(X_k(t))P_m(X_l(t))\big ) \, .
\end{align*}
Consequently by using Cauchy-Schwarz inequality and stationarity,
\begin{eqnarray*}
\frac{1}{n} \Big \| \sum_{k=1}^{n-N} [ X_k(t)-\mathbb{E}(X_k(t) \mid \mathcal{F}_n)]\Big \| ^2_2
&\leq & \frac{1}{n} \sum_{m=n+1}^{\infty} \sum_{k=1-m}^{n-N-m} \sum_{l=1-m}^{n-N-m} \| P_0(X_k(t))\| _2 \| P_0(X_l(t))\| _2 \\
&\leq & \big (\sum_{k=-\infty}^{-(N+1)} \| P_0(X_k(t))\| _2 \big )^2 \, .
\end{eqnarray*}
Hence we get
\begin{equation*}
\frac{1}{\sqrt{n}} \int_{\mathbf{T}} \Big \| \sum_{k=1}^{n-N} [X_k(t)-\mathbb{E}(X_k(t) \mid \mathcal{F}_n)] \Big \| _2 \, \mu(dt) \leq \sum_{k=-\infty}^{-(N+1)} \int_{\mathbf{T}} \| P_0(X_k(t))\| _2 \, \mu(dt),
\end{equation*}
and by $(\ref{eqnthmclt2})$, it follows that
\begin{equation}\label{eqndemthmproj11}
\underset{N\rightarrow \infty}{\lim} \ \underset{n\rightarrow\infty}{\limsup} \  \frac{1}{\sqrt{n}} \int_{\mathbf{T}} \Big \| \sum_{k=1}^{n-N}[X_k(t)-\mathbb{E}(X_k(t)\mid \mathcal{F}_n)] \Big \|_2\, \mu(dt)=0.
\end{equation}
\noindent Also by stationarity,
\begin{eqnarray*}
\Big \| \sum_{k=n-N+1}^n [X_k(t)-\mathbb{E}(X_k(t) \mid \mathcal{F}_n)] \Big \| _2 & = & \Big \| \sum_{k=n-N+1}^n \sum_{m=n+1}^{\infty} P_m(X_k(t)) \Big \| _2 \\
                                                                    & \leq & \sum_{k=n-N+1}^n \big (\sum_{m \in \mathbb{Z}} \| P_m(X_k(t))\| _2 \big ) \\
                                                                  & \leq & N \sum_{i \in \mathbb{Z}} \| P_0(X_i(t))\| _2 \, .
\end{eqnarray*}
Hence
\begin{equation} \label{eqndemthmproj12}
\int_{\mathbf{T}} \frac{1}{\sqrt{n}} \Big \| \sum_{k=n-N+1}^n [X_k(t)-\mathbb{E}(X_k(t) \mid \mathcal{F}_n)] \Big \| _2 \, \mu(dt)  \leq  \frac{N}{\sqrt{n}} \sum_{i \in \mathbb{Z}} \int_{\mathbf{T}} \| P_0(X_i(t))\| _2 \, \mu(dt).
\end{equation}
Therefore, $(\ref{eqndemthmproj11})$ and $(\ref{eqndemthmproj12})$ imply that
\begin{equation*}
\underset{n\rightarrow \infty}{\lim} \int_{\mathbf{T}} \Big\| \frac{S_n(t)}{\sqrt{n}}-\mathbb{E}\Big(\frac{S_n(t)}{\sqrt{n}} \Big| \mathcal{F}_n \Big)\Big\|_2 \, \mu(dt)=0.
\end{equation*}
To end the proof, it remains to prove (\ref{eqnthmclt4}). With this aim, we use Corollary $1$ in Dedecker, Merlev\`ede and Voln\'y \cite{dedeckermerlevedevolny07}.
Hence it suffices to prove that, for any $f$ in $\mathbf{L}^{\infty}(\mu)$,
\begin{equation}\label{eqndemthmproj13}
\sum_{k\in \mathbb{Z}} \|P_0(f(X_k))\|_2<\infty.
\end{equation}
As $f$ is a linear form on $\mathbf{L}^{1}(\mu)$ then $f$ belongs to $\mathbf{L}^{\infty}(\mu)$. It follows that
\begin{eqnarray*}
\|P_0(f(X_k))\|_2=\|f(P_0(X_k))\|_2 &\leq& C(f) \Big \| \int_{\mathbf{T}} |P_0(X_k)(t)| \, \mu(dt) \Big\|_2 \\
                                       &\leq& C(f) \int_{\mathbf{T}} \|P_0(X_k)(t)\|_2 \, \mu(dt),
\end{eqnarray*}
where $C(f)$ is a constant depending on $f$.
Consequently, $(\ref{eqndemthmproj13})$ holds as soon as $(\ref{eqnthmclt2})$ holds.
\hfill $\square$
\subsection{Proof of Corollary \ref{corclt1}}
$\newline$
\indent We prove Corollary \ref{corclt1} with the same arguments as in the end of the proof of Corollary 2 in Peligrad and Utev \cite{peligradutev06}. \\
\indent By stationarity and orthogonality, for all $k$ in $\mathbb{Z}$, we have
\begin{equation*}
\|\mathbb{E}(X_k\mid \mathcal{F}_0)\|^2_2
                                           =  \Big \|\sum_{j=-\infty}^{0} P_j(X_k) \Big \|_2^2
                                           =  \sum_{j=k}^{\infty}\|P_{-j}(X_0)\|_2^2
\end{equation*}
and
\begin{equation*}
\|X_{-k}-\mathbb{E}(X_{-k}\mid \mathcal{F}_0)\|_2^2 =  \Big \|\sum_{j=1}^{\infty} P_j(X_{-k}) \Big \|_2^2
                                                     =  \sum_{j=1}^{\infty} \|P_j(X_{-k})\|_2^2
                                                     =  \sum_{j=k+1}^{\infty} \|P_j(X_0)\|_2^2.
\end{equation*}
Now, applying Lemma $A.2$ in Peligrad and Utev \cite{peligradutev06},
to $a_i:=\|P_{-i}(X_0)\|_2$, it follows
\begin{eqnarray}
\sum_{i=1}^{\infty} \|P_{-i}(X_0)\|_2 & \leq & 3 \sum_{n=1}^{\infty}n^{-1/2}\Big(\sum_{i=n}^{\infty} \|P_{-i}(X_0)\|_2^2 \Big)^{1/2} \notag  \\
                                      & \leq & 3 \sum_{n=1}^{\infty} n^{-1/2} \|\mathbb{E}(X_n \mid \mathcal{F}_0)\|_2, \label{eqncorproof1}
\end{eqnarray}
and then to $a_i:=\|P_i(X_0)\|_2$,
\begin{eqnarray}
\sum_{i=1}^{\infty} \|P_{i}(X_0)\|_2 & \leq
  & 3 \sum_{n=1}^{\infty}n^{-1/2}\Big(\sum_{i=n}^{\infty} \|P_{i}(X_0)\|_2^2 \Big)^{1/2} \notag  \\
                                            & \leq & 3 \sum_{n=1}^{\infty} n^{-1/2} \|X_{-(n-1)}-\mathbb{E}(X_{-(n-1)}|\mathcal{F}_0)\|_2 \notag \\
                                           & \leq & 3 \sum_{n=0}^{\infty} (n+1)^{-1/2} \|X_{-n}-\mathbb{E}(X_{-n}|\mathcal{F}_0)\|_2 \notag \\
                                           & \leq & 3 \|X_0-\mathbb{E}(X_0|\mathcal{F}_0)\|_2+3 \sum_{n=1}^{\infty} n^{-1/2} \|X_{-n}-\mathbb{E}(X_{-n}|\mathcal{F}_0)\|_2 \notag  \\
                                      & \leq & 3 \sum_{n=1}^{\infty} n^{-1/2} \|X_{-n}-\mathbb{E}(X_{-n} \mid \mathcal{F}_0)\|_2 +6\|X_0\|_2. \label{eqncorproof2}
\end{eqnarray}
Therefore, by $(\ref{eqncorproof1})$ and $(\ref{eqncorproof2})$, we deduce that
\begin{eqnarray*}
 \sum_{k \in \mathbb{Z}} \|P_0(X_k)\|_2  & = & \sum_{k=1}^{\infty} \|P_{-k}(X_0)\|_2 +\sum_{k=1}^{\infty} \|P_k(X_0)\|_2+\|P_0(X_0)\|_2 \\
                                          & \leq & 3 \sum_{n=1}^{\infty} n^{-1/2} \Big( \sum_{i=n}^{\infty} \|P_{-i}(X_0)\|_2^2 \Big)^{1/2} \\
                                           &      & \qquad \qquad + 3 \sum_{n=1}^{\infty} n^{-1/2} \Big( \sum_{i=n}^{\infty}\|P_i(X_0)\|_2^2  \Big)^{1/2} + \|P_0(X_0)\|_2 \\
                                          & \leq & 3 \big[\sum_{n=1}^{\infty} n^{-1/2} \|\mathbb{E}(X_n \mid \mathcal{F}_0)\|_2 + \sum_{n=1}^{\infty}n^{-1/2}\|X_{-n}-\mathbb{E}(X_{-n}\mid \mathcal{F}_0)\|_2 \big] + 8  \|X_0\|_2.
\end{eqnarray*}
Consequently, $(\ref{eqncorclt12})$ and $(\ref{eqncorclt12bis})$ implies $(\ref{eqnthmclt2})$.
\hfill $\square$
\subsection{Proof of Proposition \ref{propphimelange}.}
$\newline$
\indent We  apply Corollary \ref{corclt1} to the variables
$X_i(.)=\{t\mapsto \mathbf{1}_{Y_i \leq t}-F_Y(t), t \in \mathbb{R} \}$. \\
Let
\begin{equation*}
Z(t)= \frac{F_{Y_k \mid \mathcal{F}_0}(t)-F_Y(t)}{\|F_{Y_k \mid \mathcal{F}_0}(t)-F_Y(t)\|_2}.
\end{equation*}
Obviously,
\begin{equation}\label{eqnalphamelange2}
\forall \ t \in \mathbb{R},  \  \|\mathbb{E}(X_k(t)\mid \mathcal{F}_0)\|_2 \leq \|X_k(t)\|_2 \leq \sqrt{F_Y(t)(1-F_Y(t))}.
\end{equation}
By Proposition $2.1$ in Dedecker \cite{dedecker04} and by (\ref{eqnalphamelange2}), for any $k$ in $\mathbb{Z}$, we derive
\begin{eqnarray*}
\|F_{Y_k \mid \mathcal{F}_0}(t)-F_Y(t)\|_2 & = & \mathbb{E} \Big ( \frac{F_{Y_k \mid \mathcal{F}_0}(t)-F_Y(t)}{\|F_{Y_k \mid \mathcal{F}_0}(t)-F_Y(t)\|_2} (F_{Y_k \mid \mathcal{F}_0}(t)-F_Y(t)) \Big) \\
                                                             & = & \mbox{Cov}(Z(t),X_k(t))\\
                                                             & \leq & 2  \|Z(t)\|_2 \|X_k(t)\|_2 \sqrt{\tilde{\phi}(k)} \\
                                                             & \leq & 2  \sqrt{F_Y(t)(1-F_Y(t))} \sqrt{\tilde{\phi}(k)}.
\end{eqnarray*}
Consequently, we deduce that (\ref{eqncorclt12}) holds as soon as (\ref{eqnpropphimelange1}) holds.
\hfill $\square$
\subsection{Proof of Proposition \ref{propalphamelange}.}
$\newline$
\indent We apply Corollary \ref{corclt1} to the random variables $X_i(.)=\{t\mapsto \mathbf{1}_{Y_i \leq t}-F_Y(t), t \in \mathbb{R} \}$. \\
By H\"older's inequality for any $k \geq 0$, we get
\begin{eqnarray}
 \|\mathbb{E}(X_k \mid \mathcal{F}_0)\|_2 & = & \|\mathbb{E}(\mathbf{1}_{Y_k \leq t} \mid \mathcal{F}_0)-\mathbb{P}(Y_k \leq t)\|_2 \notag \\
                                                          & \leq & \sqrt{\|\mathbb{E}(\mathbf{1}_{Y_k \leq t} \mid \mathcal{F}_0)-\mathbb{P}(Y_k \leq t)\|_1}\sqrt{\|\mathbb{E}(\mathbf{1}_{Y_k \leq t} \mid \mathcal{F}_0)-\mathbb{P}(Y_k\leq t)\|_{\infty}} \notag \\
                                                          & \leq &  \sqrt{\tilde{\alpha}(k)} \label{eqnalphamelange1}.
\end{eqnarray}
Using (\ref{eqnalphamelange2}) and (\ref{eqnalphamelange1}),
\begin{eqnarray*}
\int_{-\infty}^{\infty} \| \mathbb{E}(X_k (t)\mid \mathcal{F}_0)\|_2 \, dt &\leq & \int_{\mathbb{R}} \sqrt{\tilde{\alpha}(k)}\wedge \sqrt{F_{Y}(t)(1-F_Y(t))} \, dt \\
                                                                           & \leq & \int_0^{\infty} \sqrt{\tilde{\alpha}(k)} \wedge \sqrt{1-F_Y(t)} \, dt + \int_{-\infty}^0 \sqrt{\tilde{\alpha}(k)} \wedge \sqrt{F_Y(t)} \, dt \\
                                                                           & \leq & \int_0^{\infty} \sqrt{\tilde{\alpha}(k)} \wedge \sqrt{\mathbb{P}(|Y|>t)} \, dt.
\end{eqnarray*}
Notice that
\begin{eqnarray*}
\sqrt{\tilde{\alpha}(k)} \wedge \sqrt{\mathbb{P}(|Y|>t)} & = & \int_0^1 \mathbf{1}_{u \leq \sqrt{\mathbb{P}(|Y|>t)}} \mathbf{1}_{u \leq \sqrt{\tilde{\alpha}(k)}} \, du \\
                                                                  & = & \int_0^{\sqrt{\tilde{\alpha}(k)}} \mathbf{1}_{u^2 \leq \mathbb{P}(|Y|>t)} \, du \\
                                                                  & = & \int_0^{\sqrt{\tilde{\alpha}(k)}} \mathbf{1}_{Q_Y(u^2) \geq t} \, du.
\end{eqnarray*}
By Fubini and by a change of variable, we derive
\begin{eqnarray*}
\int_0^{+\infty} \sqrt{\tilde{\alpha}(k)} \wedge \sqrt{\mathbb{P}(|Y|>t)} \, dt & = & \int_0^{\infty} \Big( \int_0^1 \mathbf{1}_{u \leq \sqrt{\mathbb{P}(|Y|>t)}} \mathbf{1}_{u \leq \sqrt{\tilde{\alpha}(k)}} \, du \Big) \, dt \\
                                                                                         & = & \int_0^{\sqrt{\tilde{\alpha}(k)}} \Big( \int_0^{\infty} \mathbf{1}_{Q_Y(u^2) \geq t} \, dt \Big) \, du \\
                                                                                         & = & \int_0^{\sqrt{\tilde{\alpha}(k)}} Q_Y(u^2) \, du
                                                                                          =  \frac{1}{2}\int_0^{\tilde{\alpha}(k)} \frac{Q_Y(u)}{\sqrt{u}} \, du.
\end{eqnarray*}
We deduce that (\ref{eqnpropalphamelange1}) implies (\ref{eqncorclt12}). The CLT holds, by applying Corollary \ref{corclt1}.
\hfill $\square$
\subsection{Proof of Corollary \ref{exphi}.}
$\newline$
\indent Recall that to prove the convergence for the empirical distribution function of $Y$ where $(Y_k)_{k \in {\mathbb Z}}$ is defined by (\ref{defYk}), it suffices to show the convergence in ${\mathbf L}^1(\lambda)$ of the empirical distribution function of $f(Z)$ where $(Z_i)_{i \in \mathbb{Z}}$ is the stationary Markov chain with transition Kernel $K$. Hence we shall prove that
$f(Z)$ satisfies the conditions of Proposition  \ref{propphimelange}.

Using the fact that $f$ is a monotonic function, Dedecker and Prieur \cite{dedeckerprieur05} proved that
\begin{equation*}
\tilde{\phi}(\mathcal{F}_0, f(Z_k)) \leq \tilde{\phi}(\mathcal{F}_0,Z_k).
\end{equation*}
In addition, they proved that if $(\ref{eqn110000})$ holds then
\begin{equation*}
\tilde{\phi}(\mathcal{F}_0,Z_k)\leq C_1 \rho^k,
\end{equation*}
with $C_1$ a positive constant (see \cite{dedeckerprieur05}). This entails that  $(\ref{eqnpropphimelange1})$ holds.
\hfill $\square$
\subsection{Proof of Corollary \ref{coralpha}.}
$\newline$
\indent For the same reasons given in the proof of Corollary \ref{exphi}, we shall prove that
$f(Z)$ satisfies the conditions of Proposition  \ref{propalphamelange}. Hence we apply Proposition \ref{propalphamelange} to the variables $X_i(.)=\{t \mapsto \mathbf{1}_{Y_i \leq t}-F_Y(t), t \in \mathbb{R}\}$. Using the fact that $f$ is a monotonic function, Dedecker and Prieur \cite{dedeckerprieur05} proved that
\begin{equation*}
\tilde{\alpha}(\mathcal{F}_0,Y_k)=\tilde{\alpha}(\mathcal{F}_0, f(Z_k)) \leq \tilde{\alpha}(\mathcal{F}_0,Z_k).
\end{equation*}
Recently, Dedecker, Gou\"ezel and Merlev\`ede \cite{dedeckergouezelmerlevede08} proved in Proposition $1.12$, that there exists a constant $C_{\gamma}$, such that, for any positive integer $k$,
\begin{equation} \label{eqnalphamelange3}
\tilde{\alpha}(\mathcal{F}_0,Z_k) \leq \frac{C_{\gamma}}{(n+1)^{\frac{1-\gamma}{\gamma}}}.
\end{equation}

As $T_{\gamma}$ is an intermittent map, and $f$ is a monotonic function, it follows by (\ref{eqnalphamelange3}) that
\begin{equation*}
\int_{0}^{\infty} \sqrt{\tilde{\alpha}(k)} \wedge \sqrt{\nu_{\gamma}(|f|>t)}\, dt \leq \int_{0}^{\infty} \frac{C_{\gamma}}{(k+1)^{\frac{1-\gamma}{2\gamma}}} \wedge \sqrt{\nu_{\gamma}(|f|>t)}\, dt.
\end{equation*}
Consequently, (\ref{eqnpropalphamelange2}) holds as soon as (\ref{eqnalphamixing1}) holds.
\hfill $\square$
\subsection{Proof of Remark \ref{remG}} \label{proofremG}
$\newline$
\indent Since the density $g_{\nu_{\gamma}}$ of $\nu_{\gamma}$ is such that $g_{\nu_{\gamma}}(x) \leq V(\gamma) x^{-\gamma}$, we infer that
\begin{equation} \label{eqnalphamelange4}
\nu_{\gamma}(f>t) \leq \frac{D^{\frac{1-\gamma}{a}}V(\gamma)}{1-\gamma} t^{-\frac{1-\gamma}{a}} \, ,
\end{equation}
where  $D$ and $V(\gamma)$ are positive constants.
By Fubini and
(\ref{eqnalphamelange4}), we then  get that
\begin{eqnarray*}
\int_0^{+\infty} k^{-\frac{1-\gamma}{2 \gamma}} \wedge \sqrt{\nu_{\gamma}(|f| >t)} \, dt & = & \int_0^{\infty} \Big( \int_0^1 \mathbf{1}_{u \leq \sqrt{\nu_{\gamma}(|f| >t)}} \mathbf{1}_{u \leq k^{-\frac{1-\gamma}{2\gamma}}} \, du \Big) \, dt \\
 & \leq & K(\gamma, a) k^{ \big (- \frac{(1-\gamma)}{2 \gamma} \big) \big(- \frac{2a}{1-\gamma }+ 1 \big)},
\end{eqnarray*}
where $K(\gamma,a)$ is a constant. \\
Consequently, (\ref{eqnalphamixing1}) holds as soon as $a<\frac{1}{2}-\gamma$ does.
\hfill $\square$

\subsection{Proof of Proposition \ref{propcausallinearprocessadapt}.}
$\newline$
\indent We apply Theorem \ref{thmclt} to the random variables $X_k(t)=\{t\mapsto \mathbf{1}_{Y_k\leq t }-F_Y(t)\}$. Let $\mathcal{M}_i=\sigma(\varepsilon_k,k\leq i)$.
By a result in Lemma $6$ in Dedecker and Merlev\`ede \cite{dedeckermerlevedeedf07},
\begin{equation}\label{eqnpreuvepropcausallinearprocessadapt1}
\|F_{Y_k \mid \mathcal{M}_0}(t)-F_{Y_k \mid \mathcal{M}_{-1}}(t)\|_2 \leq K \mid a_0\mid^{-1}\mid a_k\mid \|\varepsilon_1-\varepsilon_0\|_2.
\end{equation}
Moreover, we have
\begin{eqnarray*}
\|F_{Y_k \mid \mathcal{M}_0}(t)-F_{Y_k \mid \mathcal{M}_{-1}}(t)\|_2 & \leq & \|F_{Y_k \mid \mathcal{M}_0}(t)-F_Y(t)\|_2 + \|F_{Y_k \mid \mathcal{M}_{-1}}(t)-F_Y(t)\|_2 \\
                                                                                       & \leq & 2 \sqrt{F_Y(t)(1-F_Y(t))}.
\end{eqnarray*}
We deduce that (\ref{eqnthmclt2}) holds as soon as
\begin{equation*}\label{eqncausallinearprocessadapt1}
\sum_{k \in \mathbb{Z}} \int_0^{\infty} (K \mid a_0 \mid^{-1} \mid a_k \mid \|\varepsilon_1-\varepsilon_0\|_2) \wedge ( 2 \sqrt{\mathbb{P}(|Y_k|>t)}) \, dt < \infty,
\end{equation*}
and it may be reduced to
\begin{equation*}\label{eqncausallinearprocessadapt2}
\sum_{k \in \mathbb{Z}} \int_0^{\infty} \mid a_k\mid \wedge\sqrt{\mathbb{P}(|Y_k|>t)} \, dt < \infty \Leftrightarrow \sum_{k \in \mathbb{Z}} \int_0^{\mid a_k \mid ^2} \frac{Q_{\mid Y_k\mid}(u)}{\sqrt{u}} \, du < \infty.
\end{equation*}
Now, from Theorem \ref{thmclt},
we infer that $\sqrt{n}(F_n-F_Y)$ converges in law to a $\mathbf{L}^1(\lambda)$-valued centered Gaussian random variable $G$, with covariance  operator $\Phi_{\mu}$ defined by $(\ref{eqnpropphimelange10})$.
\hfill $\square$
\subsection{Proof of Remark \ref{remRio}.}
$\newline$
\indent By using Lemma $2.1$ in Rio \cite{rio00}, page $35$, we have that
\begin{eqnarray*}
\int_0^{\mid a_k \mid^2} \frac{Q_{\mid Y_0 \mid }(u)}{\sqrt{u}} \, du & \leq & \big( \sum_{j \geq 0} |a_j| \big) \int_0^{\mid a_k \mid^2} \frac{Q_{\mid \varepsilon_0 \mid}(u)}{\sqrt{u}} \, du.
\end{eqnarray*}
Consequently since $\sum_{j \geq 0} |a_j| < \infty$,  (\ref{eqnpropcausallinearprocessadapt1}) is true provided that
\begin{equation}\label{eqncausallinearprocess3bis}
\sum_{k\in \mathbb{Z}} \int_0^{\mid a_k \mid^2} \frac{Q_{\mid \varepsilon_0 \mid }(u)}{\sqrt{u}} \, du < \infty.
\end{equation}
\hfill $\square$
 \subsection{Proof of Corollary \ref{propmom}.}

\subsubsection{Proof of Item $1$ of Corollary \ref{propmom}.}
$\newline$
\indent To apply Corollary \ref{propcausallinearprocessadapt},  it suffices to prove (\ref{eqncausallinearprocess3}). \\
\indent Firstly, recall that, if $U$ is an uniform random variable on $[0,1]$, $Q^2_{|\varepsilon_0|}(U)$ and $|\varepsilon_0|^2$ have the same law. \\
We proceed as in Rio \cite{rio00} p 15.
By H\"older's inequality on $[0,1]$, it follows that
\begin{eqnarray*}
\sum_{k \geq 0} \int_0^{|a_k|^2} \frac{Q_{|\varepsilon_0|}(u)}{\sqrt{u}} \, du & = & \int_0^1 Q_{|\varepsilon_0|}(u) \Big( \frac{\sum_{k \geq 0} \mathbf{1}_{\{u \leq |a_k|^2\}}}{\sqrt{u}} \Big) \\
                                                                               & \leq & \Big(\int_0^1 Q_{|\varepsilon_0|}(u)^r \, du \Big)^{1/r} \Big(\int_0^1 \Big(\frac{\sum_{k \geq 0} \mathbf{1}_{\{u\leq |a_k|^2\}}}{\sqrt{u}}\Big)^{r/(r-1)} \, du \Big)^{(r-1)/r}.
\end{eqnarray*}
Using the same notations as in Dedecker and Doukhan \cite{dedeckerdoukhan03}, let
\begin{equation*}
\delta^{-1}(u)=\sum_{k \geq 0}\mathbf{1}_{\{u \leq |a_k|^2 \}} \ \mbox{and} \ f(x)=x^{r/(r-1)}.
\end{equation*}
We infer that
\begin{eqnarray*}
f(\delta^{-1}) & =& \sum_{j=0}^{\infty} (f(j+1)-f(j))\mathbf{1}_{\{u\leq |a_j|^2\}} \\
               & = & \sum_{j=0}^{\infty} ((j+1)^{r/(r-1)}-j^{r/(r-1)})\mathbf{1}_{\{u \leq |a_j|^2\}}.
\end{eqnarray*}
Set $C_r=1\vee(\frac{r}{r-1})$. Since $(j+1)^{r/(r-1)}-j^{r/(r-1)} \leq C_r j^{1/(r-1)}$,
$(\ref{eqncausallinearprocess3})$ holds as soon as
\begin{equation*}\label{eqnmom6}
\int_{0}^1 \sum_{j \geq 0} j^{1/(r-1)} \frac{\mathbf{1}_{\{ u \leq |a_j|^2\}}}{u^{\frac{r}{2(r-1)}}} \, du <\infty \, ,
\end{equation*}
which is true provided that
\begin{equation*}
\sum_{j \geq 0} j^{1/(r-1)} |a_j|^{\frac{r-2}{r-1}}<\infty.
\end{equation*}
\hfill $\square$

\subsubsection{Proof of Item $2$ of Corollary \ref{propmom}.}
$\newline$
We apply Corollary \ref{propcausallinearprocessadapt}, so it suffices to prove $(\ref{eqncausallinearprocess3})$. \\
Notice that, the quantile function $Q_{|\varepsilon_0|}$, here, is dominated by $cu^{-1/r}$. Thus, we derive
\begin{eqnarray*}
\int_0^{|a_k|^2} \frac{Q_{|\varepsilon_0|}(u)}{\sqrt{u}} \, du & \leq & \int_0^{|a_k|^2} \frac{c}{u^{1/2+1/r}} \, du \\
                                                              & \leq & c \frac{2r}{r-2} |a_k |^{1-2/r}.
\end{eqnarray*}
Consequently, $(\ref{eqncausallinearprocess3})$ holds as soon as $(\ref{eqnqueue})$ does.
\hfill $\square$
\nobreak
\renewcommand{\refname}{References}


\begin{thebibliography}{99}
\bibitem{deacostaaraujogine78} de Acosta, A., Araujo, A. and Gin\'e, E. (1978). On poisson meausures, gaussian measures and the central limit theroem in Banach spaces. \textit{Probability on Banach spaces. Adv. Probab. Related Topics.} \textbf{4}, 1-68.
\bibitem{billingsley61} Billingsley, P. (1961). The Lindeberg-L\'evy theorem for martingales. \textit{Proc. Amer. Math. Soc} \textbf{12}, No. 5, 788-792.
\bibitem{delbarrioginematran99} Del Barrio, E., Gin\'e, E. and Matr\'an, C. (1999). Central limit theorems for the Wasserstein distance between the empirical and the true distributions. \textit{The Annals of Probability} \textbf{27}, No. 2, 1009-1071.
\bibitem{dedecker04} Dedecker, J. (2004). In\'egalit\'es de covariance. \textit{C. R. Math. Acad. Sci. Paris.} \textbf{339}, No. 7, 503-506.
\bibitem{dedeckerdoukhan03} Dedecker, J. and Doukhan, P. (2003). A new covariance inequality and applications. \textit{Stochastic Processes and their Applications} \textbf{106}, 63-80.
\bibitem{dedeckergouezelmerlevede08} Dedecker, J., Gou\"ezel, S. and Merlev\`ede, F. (2008). Some almost sure results for unbounded functions of intermittent maps and their associated Markov chains.
\bibitem{dedeckermerlevedeedf07} Dedecker, J. and Merlev\`ede, F. (2007). The empirical distribution function for dependent variables: asymptotic and non asymptotic results in $\mathbf{L}^p$. \textit{ESAIM Probability and Statistics.} \textbf{11}, 102-114.
\bibitem{dedeckermerlevedevolny07} Dedecker, J., Merlev\`ede, F. and Voln\'y, D. (2007). On the weak invariance principle for non-adapted sequences under projective criteria. \textit{J. Theoret. Probab.} \textbf{20}, No. 4, 971-1004.
\bibitem{dedeckerprieur05} Dedecker, J. and Prieur, C. (2005). New dependence coefficients. Examples and applications to statistics. \textit{Probab. Theory. Relat. Fields} \textbf{132}, 203-236.
\bibitem{gouezel04} Gou\"ezel, S. (2004). Central limit theorem and stable laws for intermittent maps. \textit{Probab. Theory Relat. Fields} \textbf{128}, 82-122.
\bibitem{hennionherve01} Hennion, H. and Herv\'e, L. (2001). Limit theorems for Markov chains and stochastic properties of dynamical systems by quasi-compactness. \textit{Lecture Notes in Mathematics.} \textbf{1766}.
\bibitem{ibragimov62} Ibragimov, I.A. (1962). Some limit theorems for stationary processes. \textit{Theory Probab. Appl.} \textbf{7}, 349-382.
\bibitem{jain77} Jain, N. (1977). Central limit theorem and related questions in Banach space. \textit{Proceeding of symposia in Pure Mathematics XXXI.} 55-65.
\bibitem{liveranisaussolvaienti99} Liverani, C., Saussol, B. and Vaienti, S. (1999). A probabilistic approach to intermittency. \textit{Ergodic Theory Dynam. Systems} \textbf{19}, 671-685.
\bibitem{peligradutev06} Peligrad, M. and Utev, S. (2006). Central limit theorem for stationary linear processes. \textit{Annals of Probability} \textbf{34}, No. 4, 1608-1622.
\bibitem{pomeaumanneville80} Pomeau, Y., Manneville, P. (1980). Intermittent transition to turbulence in dissipative dynamical systems. \textit{Commun. Math. Phys.} \textbf{74}, 189-197.
\bibitem{rio00} Rio, E. (2000). Th\'eorie asymptotique des processus al\'eatoires faiblement d\'ependants. \textit{Math\'ematiques et applications de la SMAI.} \textbf{31}.
\bibitem{rosenblatt56} Rosenblatt, M. (1956). A central limit Theorem and a strong mixing condition. \textit{Proc. Nat. Ac. SC. U.S.A.} \textbf{42}, 43-47.
\bibitem{volny93} Voln\'y, D. (1993). Approximating martingales and the central limit theorem for strictly stationary processes. \textit{Stochastic Processes and their Applications.} \textbf{44}, 41-74.
\end{thebibliography}
\end{document}